\documentclass[12pt]{article}
\usepackage{amsmath}
\usepackage{amsfonts,amssymb}
\usepackage{epsfig}
\usepackage{graphicx}
\newcommand{\R}{{\mathbb R}}

\newcommand{\B}{\bigskip}
\newcommand{\m}{\medskip}
\newcommand{\proof}{\noindent{\em Proof: }}

\newcommand{\qed}{\hspace{\fill}$\square$}
\newtheorem{theorem}{Theorem}
\newtheorem{cor}{Corollary}
\newtheorem{lemma}{Lemma}

\newtheorem{example}{Example}
\newtheorem{definition}{Definition}

\def\sqr#1#2{{\vcenter{\vbox{\hrule height.#2pt
 \hbox{\vrule width.#2pt height#1pt \kern#1pt
 \vrule width.#2pt}
 \hrule height.#2pt}}}}
\def\square{\mathchoice\sqr68\sqr68\sqr{2.1}3\sqr{1.5}3}
\date{ }

\title{The Eigenvalue Problem for Linear and Affine Iterated Function Systems} 
\author{Michael Barnsley$^a$ \\$^a$  Mathematical Sciences Institute,
Australian National University \\ Canberra, ACT 0200, Australia \\ {\tt Michael.Barnsley@anu.edu.au} \\ \\
Andrew Vince$^b$ (corresponding author) \\  $^b$ University of Florida, Department of Mathematics \\ 358 Little Hall, PO Box 118105, Gainesville, FL 32611-8105, USA \\ {\tt  avince@ufl.edu} \\  phone:  352-392-0281 ext 246, fax: 352-392-4958 }

\begin{document}
\maketitle
 
\begin{abstract} The eigenvalue problem for a linear function $L$ 
centers on solving the eigen-equation $Lx = \lambda \, x$.  This paper
generalizes the eigenvalue problem from a single linear function to an
iterated function system $F$ consisting of possibly an infinite number
of linear or affine  functions. The eigen-equation becomes 
$F(X) = \lambda \, X$,
where $\lambda >0$ is real, $X$ is a compact set, and $F(X) =
\bigcup_{f\in F} f(X)$.  The main result is that an irreducible,
linear iterated function system $F$ has a unique eigenvalue $\lambda$
equal to the joint spectral radius of the functions in $F$ and a
corresponding eigenset $S$ that is centrally symmetric, star-shaped,
and full dimensional. Results of Barabanov and of
Dranishnikov-Konyagin-Protasov on the joint spectral radius follow as
corollaries.
\end{abstract} 


{\bf Mathematical subject codes}:  15A18; 28A80

\section{Introduction}

Let $L : \R^2 \rightarrow \R^2$ be a linear map with no nontrivial
invariant subspace, equivalently no real eigenvalue.  We use the
notation $L(X) := \{Lx \, : \, x \in X\}$. Although $L$ has no real
eigenvalue, $L$ does have an eigen-ellipse.  By {\it
eigen-ellipse} we mean an ellipse $E$, centered at the origin, such
that $L(E) = \lambda\,E$, for some real $\lambda>0$.  An example of an
eigen-ellipse appears in Example~\ref{ex0} of Section~\ref{examples}  
and in Figure~\ref{f0}. 
Although easy to prove,
the existence of an eigen-ellipse appears not to be well known.

\begin{theorem} If $L : \R^2 \rightarrow \R^2$ is a linear map 
with no real eigenvalue, then there is an ellipse
$E$ and a $\lambda >0$ such that $L(E) = \lambda\,E$.
\end{theorem}

\proof Using the real Jordan canonical form for $L$, there exists
an invertible $2\times 2$ matrix $S$ such that
$$M := S^{-1} L S = \lambda \, \left ( \begin{array}{rr} \cos \, \theta &
  -\sin \, \theta \\ \sin \, \theta & \cos \,\theta \end{array} \right
),$$ for some angle $\theta$ and $\lambda >0$. If $D$ is the unit
disk centered at the origin and if $E = S(D)$, then
$$L(E) = SMS^{-1}(E) = SM(D) = \lambda S(D) = \lambda E.$$ 
\qed

The intent of this paper is to investigate the existence of eigenvalues
and corresponding eigensets in a more general setting.

\begin{definition}[iterated function system]
Let $\mathbb{X}$ be a complete metric space. If $f_i:\mathbb{X}
\rightarrow\mathbb{X}, \, i \in I$, are continuous mappings, then
$F=\left( \mathbb{X}; f_i,\, i\in I \right)$ is called an
\textbf{iterated function system} (IFS).  The set $I$ is the
index set.  Call IFS $F$ \textbf{linear}
if ${\mathbb X} = \R^n$ and each $f \in F$ is a linear map and
\textbf{affine} if  ${\mathbb X} = \R^n$ and each $f \in F$ is an affine map.
\end{definition}

In the literature the index set $I$ is usually finite.  This is
because, in constructing deterministic fractals, it is not practical
to use an infinite set of functions.  We will, however, allow an
infinite set of functions in order to obtain certain results on the
joint spectral radius. In the case of an infinite linear IFS $F$ we
will always assume that the set of functions in $F$ is compact.  For
linear maps, this just means, regarding each linear map as an $n\times
n$ matrix, that the set $F$ of linear maps is a compact subset of
$\R^{n\times n}$. \m

Let $\mathbb H = \mathbb H({\mathbb X})$ denote the collection of all 
nonempty compact subsets
of ${\mathbb X}$, and, by slightly abusing the notation, let $F \,:\,
\mathbb H({\mathbb X}) \rightarrow \mathbb H({\mathbb X})$ also denote
the function defined by
$$F(B) = \bigcup_{f\in F}\, f(B).$$  Note that, if $B$ is compact
and $F$ is compact, then $F(B)$ is also compact.
Let $F^{k}$ denote $F$
iterated $k$ times with $F^0(B) = B$ for all $B$. 
Our intention is to investigate solutions to the eigen-equation

\begin{equation} \label{eigen-equation} F(X) = \lambda \, X, \end{equation}
where $\lambda \in \R, \; \lambda >0$, and $X \neq \{0\}$ is a compact
set in Euclidean space.

\begin{definition}[eigenvalue-eigenset] The value $\lambda$ in
Equation (1) above will be called an \textbf {eigenvalue} of $F$, and
$X$ a corresponding \textbf {eigenset}. \end{definition}

When $F$ consists of a single linear map on $\R^2$, the eigen-ellipse
is an example of an eigenset.  Section~\ref{examples} contains other
examples of eigenvalues and eigensets of linear IFSs.
Section~\ref{contractiveIFS} contains background results on the joint
spectral radius of a set of linear maps and on contractive IFSs. Both
of these topics are germane to the investigation of the IFS eigenvalue
problem.  Section~\ref{linear} contains the main result on the
eigenvalue problem for a linear IFS.  

\begin{theorem} \label{T-eigen} A compact, irreducible, linear 
IFS $F$ has exactly one eigenvalue which is equal to the
joint spectral radius $\rho(F)$ of $F$.  There is a
corresponding eigenset that is centrally symmetric, star-shaped, and
full dimensional.
\end{theorem}

If $F = \{\R^n ; \, f_i,\, i\in I \}$ is an IFS, let $F_{\lambda} :=
\frac{1}{\lambda} \, F = \{\R^n ; \, \frac{1}{\lambda} f_i,\, i\in I
\}$.  Another way to view the above theorem is to consider the family
$\{ F_{\lambda} \, : \, \lambda > 0 \}$ of IFSs.  If $\lambda >
\rho(F)$, then the attractor of $F_{\lambda}$ (defined formally in
the next section) is the trivial set $\{0\}$.  If $\lambda <
\rho(F)$, then $F_{\lambda}$ has no attractor. So $\lambda =
\rho (F)$ can be considered as a ``phase transition'', at which
point a somewhat surprising phenomenon occurs - the emergence of the
centrally symmetric, star-shaped eigenset. \m

Theorems of Dranisnikov-Konyagin-Protasov and of Barabanov follow as
corollaries of Theorem~\ref{T-eigen}.  These results are discussed in
Section~\ref{Cor}.  \m

No such transition phenomenon occurs in the case of an affine, but not
linear, IFS.  A result for the affine case is the following, whose
proof appears in Section~\ref{affSec}.  

\begin{theorem} \label{affThm} For a compact, irreducible, affine, but
not linear, IFS $F$, a real number $\lambda > 0$ is an eigenvalue if
$\lambda > \rho(F)$ and is not an eigenvalue if $\lambda < \rho(F)$.
There are examples 
where $\rho(F)$ is an eigenvalue and examples 
where it is not.
\end{theorem}

The transition phenomenon resurfaces in the context of
projective IFSs, which will be the subject of a subsequent
paper.

\section{Examples} \label{examples}

\begin{figure}
\begin{center}
\vskip -35mm
\includegraphics[width=10cm, keepaspectratio]{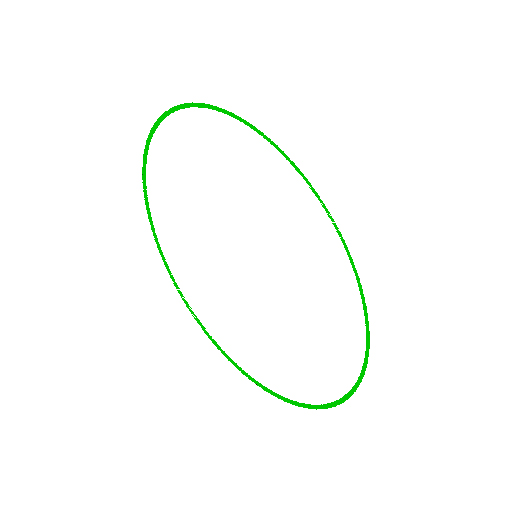}
\vskip -12mm
\caption{The eigen-ellipse for Example~\ref{ex0}}
\label{f0}
\end{center}
\end{figure}

\begin{example} \label{ex0}  Figure~\ref{f0} shows the
eigen-ellipse for the the IFS $F = (\R^2; L)$, where
$$L = \begin{pmatrix} 65.264 & -86.116 \\ 156.98 & 62.224 
 \end{pmatrix}.$$ The eigenvalue is approximately $97.23$.
\end{example}

\begin{figure}
\begin{center}
\vskip -25mm
\includegraphics[width=14cm, keepaspectratio]{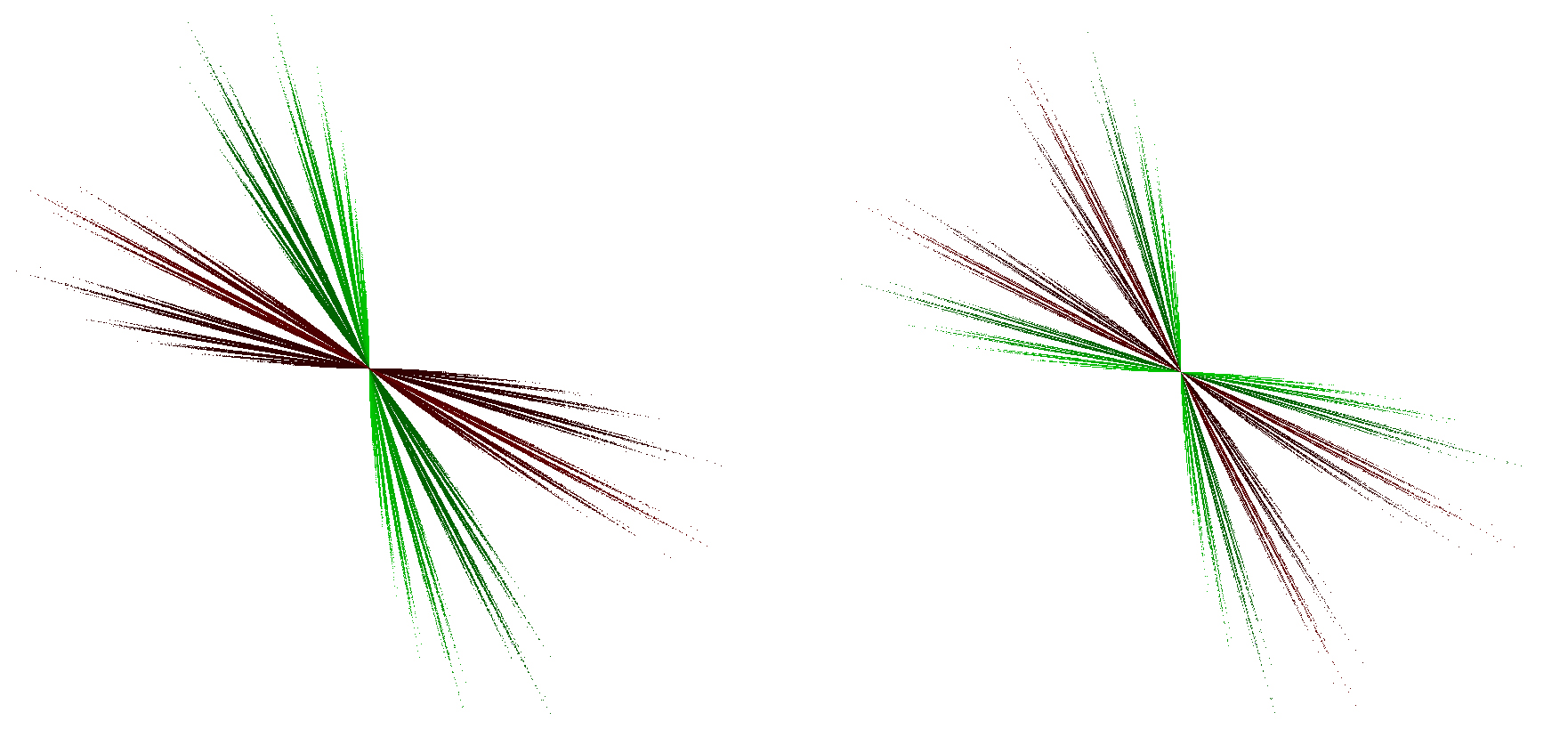}
\caption{The eigenset of Example~\ref{ex1}}
 \label{f1}
\end{center}
\end{figure}

\begin{example} \label{ex1} Figure~\ref{f1} shows an
eigenset for the IFS $F = (\R^2; L_1, L_2)$, where
$$L_1 = \begin{pmatrix} 10 & 10 \\ 8 & 0 \end{pmatrix}, \qquad \qquad
L_2= \begin{pmatrix} 8 & 0 \\ 10 & 10 \end{pmatrix}.$$ The eigenvalue
is approximately $14.9$. The picture on the right is the image of the
picture of the left after applying both transformations, then shrinking
the result about its center by a factor $14.9$. The green and brown
colors help to show how the image is acted on by the two
transformations. The dots are an artifact of rounding errors, and
serve to emphasize that the pictures are approximate.
\end{example}

\begin{figure}
\begin{center}
\vskip -35mm
\includegraphics[width=14cm, keepaspectratio]{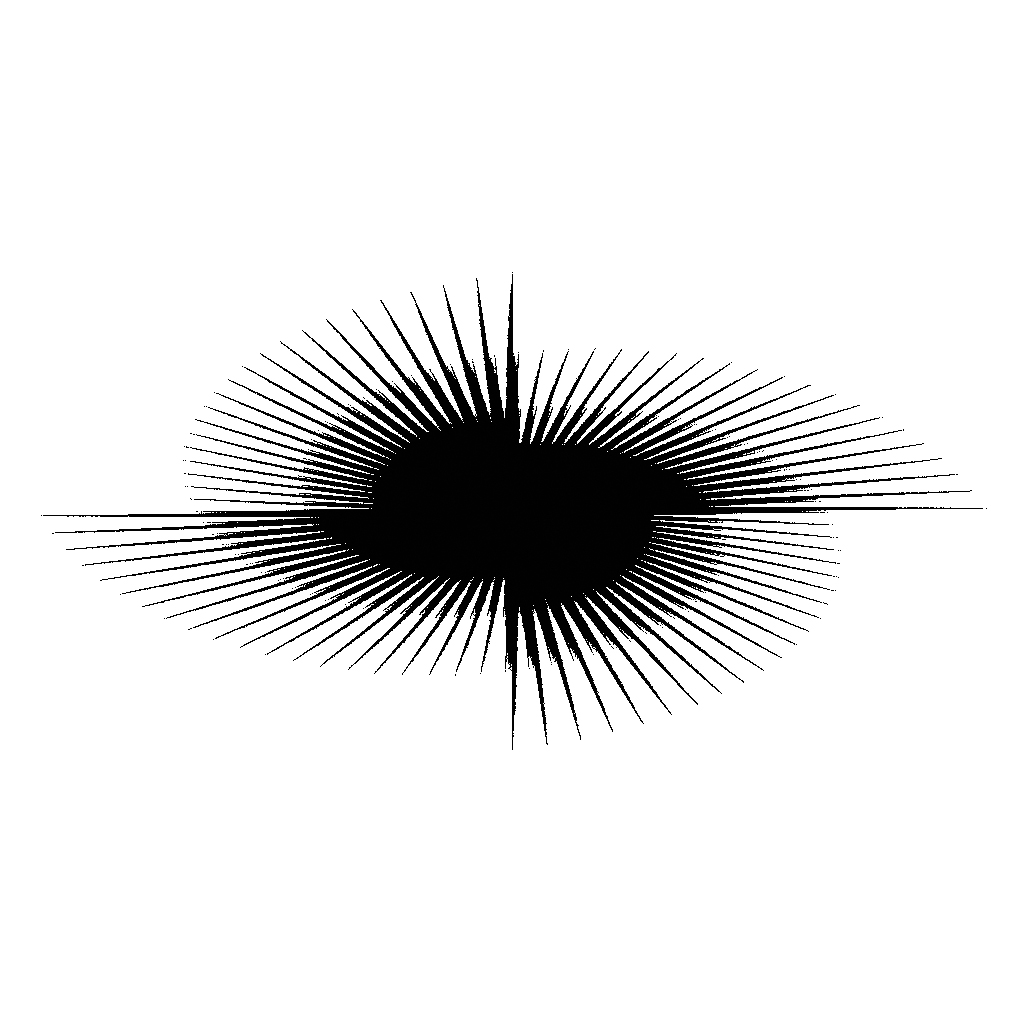}
\vskip -30mm
\caption{The eigenset of Example~\ref{ex2}}
 \label{f2} 
\end{center}
\end{figure}

\begin{example} \label{ex2}  Figure~\ref{f2} shows the
eigenset for the the IFS $F = (\R^2; L_1, L_2)$, where
$$L_1 = \begin{pmatrix}  0.02 & 0 \\0 & 1 \end{pmatrix}, \qquad \qquad
L_2= \begin{pmatrix}  0.0594 & -1.98 \\ 0.495 & 0.01547 \end{pmatrix}.$$ 
The eigenvalue is approximately $1$.
\end{example}

\section{Background} \label{contractiveIFS}

This section concerns the following three basic notions:
(1) the joint spectral radius of an IFS, (2) contractive properties of an IFS,
 and (3) the attractor of an IFS. Theorems~\ref{linearCor} and 
\ref{aff} provides the
relationship between these three notions for a linear and an affine IFS,
respectively.  

\subsection{Norms and Metrics} Any vector norm $\| \, \cdot \, \|$ on $\R^n$ 
induces a matrix norm on the space of linear maps taking $\R^n$ to
$\R^n$:
$$\| L  \| = \max \left\{ \frac{\| Lx  \|}{\| x  \|} \, :\, x \in \R^n
\right \}.$$
Since it is usually clear from the context, we use the same
notation for the vector norm as for the matrix norm. This induced 
norm is {\it sub-multiplicative}, i.e., 
$\| L\circ L'\| \leq \| L\| \cdot \| L'\|$ for any linear maps
$L,L'$. \m 

Two norms $\| \cdot \|_1$ and $\| \cdot \|_2$ are {\it equivalent} if
there are positive constants $a,b$ such that $a \| x \|_1 \leq \| x
\|_2 \leq b \| x \|_1$ for all $x\in \R^n$. Two metrics $d_1( \cdot,
\cdot)$ and $d_2( \cdot, \cdot)$ are {\it equivalent} if there exist
positive constants $a,b$ such that $a \, d_1(x,y) \leq d_2(x,y) \leq
b\, d_1(x,y)$ for all $x,y \in \R^n$.  It is well known that any two
norms on $\R^n$ are equivalent \cite{GV}. This implies that any two
$n\times n$ matrix norms are equivalent. Any norm $\| \cdot \|$ on
$\R^n$ induces a metric $d(x,y) = \| x-y \|$. Therefore any
two metrics induced from two norms are equivalent. \m 

A set $B \subset \R^n$ is called {\it centrally symmetric} if $-x\in
B$ whenever $x\in B$. A {\it convex body} in $\R^n $ is a convex set
with nonempty interior. If $C$ is a centrally symmetric convex body,
define the {\it Minkowski functional} with respect to $C$ by
$$\| x \|_C =  \inf \, \{ \mu \geq 0 \, : \, x \in \mu C\}.$$
The following result is well known.  

\begin{lemma} \label{Minkowski} The Minkowski functional is a 
norm on $\R^n$.  Conversely, any norm $\| \cdot \|$ on $\R^n$ is the
Minkowski functional with respect to the closed unit ball $\{ x \, :
\, \| x \| \leq 1\}$.
\end{lemma}

Given a metric $d(\cdot, \cdot)$, there is a corresponding 
metric $d_{\mathbb H}$, called the {\it Hausdorff metric}, on the
collection $\mathbb H(\R^n)$ of all non-empty compact subsets of
$\R^n$:
$$d_{\mathbb H}(B,C) = \max \left\{ \sup_{b\in B}\, \inf_{c\in C} d(b,c), 
\, \sup_{c\in C} \, \inf_{b\in B} d(b,c) \right\}.$$

\subsection{Joint Spectral Radius} The joint spectral radius of a 
set $\mathbb L = \{ L_i,\, i\in I\}$ of linear maps was introduced by
Rota and Strang \cite{RS} and the generalized spectral radius by
Daubechies and Lagarias \cite{DL}.  Berger and Wang \cite{BW} proved
that the two concepts coincide for bounded sets of linear maps.  The
concept has received much attention in the recent research literature;
see for example the bibliographies of \cite{S} and \cite{T}. What
follows is the definition of the joint spectral radius of $\mathbb L$.
Let $\Omega_k$ be the set of all words $i_1 \, i_2 \, \cdots \,i_k$,
of length $k$, where $i_j \in I, \, 1 \leq j \leq k$.  For $\sigma =
i_1 \, i_2 \, \cdots \, i_k \in \Omega_k$, define
$$L_{\sigma} := L_{i_1} \circ L_{i_2} \circ \cdots \circ L_{i_k}.$$ A
set of linear maps is {\it bounded} if there is an upper bound on 
their norms. Note that if $\mathbb L$ is compact, then $\mathbb L$ is
bounded.  For a linear map $L$, let $\rho (L)$ denote the ordinary
spectral radius, i.e., the maximum of the moduli of the eigenvalues of
$L$.

\begin{definition}  For any set $\mathbb L$ of linear maps and any
 sub-multiplicative norm,  the \textbf {joint spectral radius} of 
$\mathbb L$ is 
$$ \hat{\rho} = \hat{\rho}(\mathbb L) :=
\limsup_{k\rightarrow \infty} \hat{\rho}_k ^{1/k}
\qquad \mbox{where} \qquad
\hat{\rho}_k := \sup_{\sigma \in \Omega_k} \, \|L _{\sigma} \|.$$  
The \textbf {generalized spectral radius} of $\mathbb L$ is
$$\rho = \rho(\mathbb L) := \limsup_{k\rightarrow \infty}
\rho_k ^{1/k} \qquad \mbox{where} \qquad 
\rho_k := \sup_{\sigma \in \Omega_k} \, \rho(L _{\sigma}).$$
\end{definition}

\noindent The following are well known properties of the joint and
generalized spectral radius \cite{T}.
\begin{enumerate}
\item The joint spectral radius is independent of the particular
sub-multiplicative norm.
\item  For an IFS consisting of a single
linear map $L$, the generalized spectral radius is the ordinary spectral
radius of $L$.
\item For any real $\alpha >0$ we have $\rho (\alpha \, \mathbb L) =
\alpha \, \rho(\mathbb L)$ and $\hat{\rho}\, (\alpha \,  \mathbb L) = 
\alpha \,\hat{\rho}( \mathbb L)$.
\item  For all $k\geq 1$ we have 
$$\rho_k ^{1/k} \leq \rho \leq \hat{\rho} \leq \hat{\rho}_k ^{1/k},$$ 
independent of the norm used to define $\hat{\rho}$.
\item If $\mathbb L$ is bounded, then the joint and generalized
spectral radius are equal.
\end{enumerate}

From here on we always assume that $\mathbb L$ is bounded.  So, in
view of property 5, we denote by $\rho(\mathbb L)$ the common value
of the joint and generalized spectral radius.  \m

If $F$ is an affine IFS, then each $f\in F$ is of the form $f(x) = Lx + a$,
where $L$ is the {\it linear part} and $a$ is the {\it translational
part}.  Let $\mathbb L_F$ denote the set of linear parts of $F$.  

\begin{definition} The {\bf joint spectral radius of an affine IFS $F$}
is the  joint spectral radius of the set $\mathbb L_F$ of linear parts of 
$F$ and
is denoted $\rho(F)$. 
\end{definition}

\begin{definition}
A set $\{ L_i, \, i\in I\}$ of linear maps is called {\bf
reducible} if these linear maps have a common nontrivial invariant
subspace. The set is {\bf irreducible} if it is not reducible. 
An IFS is {\bf reducible (irreducible)} if the set of linear
parts is reducible (irreducible).
\end{definition}

As shown in \cite{T}, a set of linear maps is reducible if and
only if there exists an invertible matrix $T$ such that each $L_i$ can
be put simultaneously in a block upper-triangular form:
$$T^{-1} L_i T = \left ( \begin{array}{rl} A_i & * \\ 0 & B_i
\end{array} \right ),$$ with $A_i$ and $B_i$ square, and $*$ is any
matrix with suitable dimensions.  The joint spectral radius $\rho(F)$
is equal to $\max \, (\rho (\{ A_i\}), \rho(\{B_i\}))$.

\subsection{A Contractive IFS}

\begin{definition}[contractive IFS]
A function $f:\mathbb{X}\rightarrow\mathbb{X}$ is a
\textbf{contraction} with respect to a metric $d$ if there is an 
$s, \, 0\leq s<1$, such that $d(f(x),f(y))\leq s \,d(x,y)$ for all
$x,y\in{\mathbb{R}}^{n}$.  An IFS $F= ( \mathbb{X}; \, f_i,\, i\in I )$
is said to be \textbf{contractive} if there is a metric $d\,:\, \R^n
\times \R^n \rightarrow \lbrack0,\infty)$, equivalent to the standard
metric on $\R^n$, such that each $f\in F$ is a contraction 
with respect to $d$.
\end{definition}

\begin{definition}[attractor]
A nonempty compact set $A\subset\R^n$ is said to be an
\textbf{attractor} of the affine IFS $F$ if
\begin{enumerate}
\item $F(A)=A$ and
\item $\lim_{k\rightarrow\infty}F^{k}(B)=A,$ for all compact sets
$B\subset \R^n$, where the limit is with respect to the Hausdorff
metric. 
\end{enumerate}
\end{definition}

Basic to the IFS concept is the relationship between the existence of
an attractor and the contractive properties of the functions of the
IFS. The following result makes this relationship explicit in the case
of a linear IFS.  A proof of this result for an affine, but finite, IFS
appears in  \cite{ABVW}.  For completeness we provide the
proof for the infinite linear case.  The notation {\it int(X)} will be
used to denote the interior of a subset $X$ of $\R^n$.  The notation
$conv(X)$ is used for the convex hull of the set $X$. In $\R^n$ the
Minkowski sum and scalar product are defined by $Y + Z = \{y + z \, :
\, y\in Y, \, z\in Z \}$ and $\alpha \, Y = \{ \alpha \, y \, : \, y
\in Y \}$, respectively.

\begin{theorem} \label{linearCor} For a compact, linear IFS 
$F = \left( \R^n; L_i,\, i\in I \right)$
the following statements are equivalent.
\begin{enumerate}
\item {\bf [contractive]} There exists a norm $\| \cdot  \|$ on $\R^n$ 
and an $0\leq s < 1$ such that $\| Lx \| \leq s \, \|x\|$
for all $L \in F$ and all $x \in \R^n$. 
\item {\bf [F-contraction]} The map $F\,:\, \mathbb H(\R^n) \rightarrow
 \mathbb H(\R^n)$ defined by $F(B) = \bigcup_{L\in F} L(B)$ is a
 contraction with respect to a Hausdorff metric.
\item {\bf [topological contraction]} There is a compact,
 centrally symmetric, convex body $C$ such that $F(C) \subset int(C)$.
\item {\bf [attractor]} The origin is the unique attractor of $F$.
\item {\bf [JSR]} $\rho(F) < 1$.    
\end{enumerate}
\end{theorem}

\proof
(attractor $\Rightarrow$ topological contraction) Let $A$ be the
attractor of $F$.  Let $A_{\rho}=\{x\in
\mathbb{R}^{m}:d_{\mathbb{H}}(\left\{  x\right\}  ,A)\leq\rho\}$ denote the
dilation of $A$ by radius $\rho>0$. 
By the definition of the attractor, 
$\lim_{k\rightarrow\infty}d_{\mathbb{H}}(F^{k}(A_{\rho
}),A) = 0,$  so there is an integer $m$ so that 
$d_{\mathbb{H} }(F^{m}(A_{1}),A) < 1$. 
Thus, 
$$F^m (A_{1})\subset int(A_{1}).$$
\noindent If $C_1 := conv(  A_1 - A_1)$, then it is 
straightforward to check that $C_1$ is a centrally symmetric
convex body and that  $F^m (C_{1})\subset int(C_{1})$, which
implies $$conv \, F^m (C_{1})\subset int(C_{1}).$$
Consider the Minkowski sum
$$C := \sum_{k=0}^{m-1} conv \, F^k(C_{1}).$$
For any $L \in F$
\begin{align*}
L(C) 
& =
\sum_{k=0}^{m-1} L\Bigl( conv \, F^k(C_1) \Bigr)  = 
\sum_{k=0}^{m-1} conv \Bigl( L \left( F^k\left(C_{1}\right) 
\right)\Bigr) 
\\
&  \subseteq 
\sum_{k=0}^{m-1}
 conv \, F^{k+1} \left( C_{1}\right)   =
 conv \, F ^{m}(C_{1}) 
+ \sum_{k=1}^{m-1} conv \, F^k(C_{1}) \\
& \subseteq int\left( C_{1} \right) + 
\sum_{k=1}^{m-1} conv\, F^{k}(C_{1}) \\
& = 
int(C).
\end{align*}
The last equality follows from the fact that
if $K$ and $K'$ are convex bodies 
in $\mathbb{R}^{n},$ then $int(K)+ K' 
= int\left(K + K' \right)$. 
\m

(topological contraction $\Rightarrow$ contractive) Let $C$ be a
centrally symmetric, convex body such that $F(C) \subset int (C)$.
Let $\| \cdot \|_C$ be the Minkowski functional with respect to $C$
and $d_C$ the metric corresponding to the norm $\| \cdot \|_C$.  Let
$L\in F$.  Since $C$ is compact, the containment $L(C) \subset int(C)$
implies that there is an $s \in [0,1)$ such that $\|Lx \|_C \leq s \,
\|x\|_C$ for all $x\in \R^n$.  Therefore $d_C(L(x),L(y))
=\|L(x)-L(y)\|_C = \|L(x-y)\|_C \leq s\, \|x-y\|_C =s \, d_C(x,y)$,
and so $d_C$ is a metric for which each function in the IFS is a
contraction.  Since any convex body contains a ball of radius $r$ and
is contained in a ball of radius $R$ for some $r,R >0$, the metric
$d_C$ is equivalent to the standard metric. \m

(contractive $\Rightarrow$ F-contraction) In the case of an IFS $F =
(\R^n; f_i,\, i\in I)$, where $I$ is finite (and the $f_i$ are assumed
only to be continuous), this is a basic result whose proof can be
found is most texts on fractal geometry, for example \cite{F}. Since
$F$ is assumed contractive,
$$\sup \left\{ \frac{d(f_i(x),f_i(y))}{d(x,y)} \, : \, x\neq y \right\} = 
s_i < 1,$$ for each
$i \in I$.  The only sticking point in extending the proof for
the finite IFS case to the infinite IFS case is to show
that $\sup \{ s_i \, : \,i \in I \} < 1$. But if there is a sequence
$\{s_k\}$ such that $\lim_{k\rightarrow \infty} s_k = 1$, then, by the
compactness of $F$, the limit $f := \lim_{k\rightarrow \infty} f_k \in
F$. Moreover,
$$\frac{d(f(x),f(y))}{d(x,y)} = \lim_{k\rightarrow \infty}
\frac{d(f_k(x),f_k(y))}{d(x,y)} = \lim_{k\rightarrow \infty} s_k =1,$$
contradicting the assumption that each function in $F$ is a
contraction.  \m

(F-contraction $\Rightarrow$ attractor) The existence of a unique 
attractor follows directly from the Banach contraction mapping
theorem.  When $F$ is linear, uniqueness immediately
implies that the attractor is $\{0\}$. \m 

(contractive $\Leftrightarrow$ JSR) First assume that $F$ is
contractive. Hence there is an $0\leq s < 1$ such that $\| L x \| \leq
s\, \|x\|$ for all $x \in \R^n$ and all $L \in F$.  By 
property (4) of the joint spectral radius
$$\rho(F) \leq \hat{\rho_1} = \sup_{L\in F} \, 
\frac{\| L x \|}{\|x\|} \leq s <1.$$ 
The last inequality is a consequence of the compactness of $F$, the
argument identical to the one used above in showing that 
(contractive $\Rightarrow$ attractor).

Conversely, assuming
$$\limsup_{k\rightarrow \infty} \hat{\rho}_k ^{1/k} = \rho(F) < 1,$$
we will show that $F$ has attractor $A = \{0\}$.
The inequality above implies that there is an
$s$ such that $\hat {\rho_k}^{1/k} \leq s < 1$
for all but finitely many $k$.  In other words
$$\sup_{\sigma \in \Omega_k} \| L_{\sigma}\| = \hat {\rho_k} \leq s^k$$ 
for all but finitely many $k$.  For $k$ sufficiently large, this in turn
implies, for any $x \in \R^n$ and any $\sigma \in \Omega_k$, that
$\|L_{\sigma} x\| \leq s^k \|x \|$.  Therefore, for any compact set $B\subset \R^n$,
with respect to the Hausdorff metric, $\lim_{k\rightarrow \infty} F^k (B) = \{ 0 \}$.  
So $\{0\}$ is the attractor of $F$.
\qed

\begin{cor} \label{unique} If a compact, linear IFS $F$ is contractive and
$F(A) = A$ for $A$ compact, then $A = \{0\}$.
\end{cor}

\proof According to Theorem~\ref{linearCor} the IFS has the
$F$-contractive property.  According to the Banach fixed
point theorem, $F$ has a unique invariant set, i.e., a unique
compact $A$ such that $F(A) = A$. Since $F$ is linear, 
clearly $F( \{0\}) =  \{0\}$.   
\qed \m

The following theorem is an extension of Theorem~\ref{linearCor} to
the case of an affine IFS. The proof of the equivalence of the
first three statements, for a finite affine IFS, appears
in \cite{ABVW}.  The modifications in the proof (of the equivalence
of the first three statements) needed to go from the
finite to the compact case is omitted since it is exactly as in the
proof of Theorem~\ref{linearCor}.  The proof of the
equivalence of statement (4) is given below.  Note 
that this last equivalence implies that, if a linear
IFS F has an attractor and $F'$ is obtained from
$F$ by adding any translational component to each
function in $F$, then $F'$ also has an attractor
\newpage

\begin{theorem} \label{aff}
If $F = \left( \R^n; f_i,\, i\in I \right)$ is a compact, affine IFS,
then the following statements are equivalent.
\begin{enumerate}
\item {\bf [contractive]} The IFS $F$ is contractive on $\R^n$.
\item {\bf [topological contraction]} There exists a compact set $C$
such that $F(C) \subset int(C)$.
\item {\bf [attractor]} $F$ has a unique attractor, the
basin of attraction being $\R^n$.
\item {\bf [JSR]} $\rho(F) < 1$. 
\end{enumerate}
\end{theorem}

\proof  As explained above, we  prove only
 the equivalence of statement (4) to the other
statements.  Assuming  $\rho(F) < 1$ we will show that $F$ is
contractive. Let $F'$ be the linear IFS obtained from $F$
by removing the translational component from each
function in $F$.   By  Theorem~\ref{linearCor}, the IFS $F'$
is contractive.  Hence there is a norm $\| \, \cdot \, \|$ with
respect to which each $L\in F'$ is a contraction.  Define 
a metric by $d(x,y) = \| x - y \|$ for all $x,y \in \R^n$.
For any $f(x) = Lx + a \in F$ we have
$d( f(x),f(y) ) = \| f(x)-f(y)\| = \|(Lx+a) - (Ly+a)\| = \|L(x-y)\|$.
Therefore each function $f\in F$ is a contraction with respect
to metric $d$.

Conversely, assume that the affine IFS F is contractive.
With linear IFS $F'$ as defined above, 
it is shown in \cite[Theorem 6.7]{ABVW} that there
is a norm with respect to which each $L\in F'$ is
a contraction.  It follows from  Theorem~\ref{linearCor} 
that  $\rho(F) < 1$. 
\qed

\section{The Eigenvalue Problem for a Linear IFS} \label{linear}

Just as for eigenvectors of a single linear map, an eigenset of an IFS
is defined only up to scalar multiple, i.e., if $X$ is an eigenset,
then so is $\alpha \, X$ for any $\alpha >0$.  Moreover, if $X$ and
$X'$ are eigensets corresponding to the same eigenvalue, then $X\cup
X'$ is also a corresponding eigenset. For an eigenvalue of a linear
IFS, call a corresponding eigenset $X$ {\it decomposable} if $X = X_1
\cup X_2$, where $X_1 \neq X$ and $X_2 \neq X$ are also corresponding
eigensets.  Call eigenset $X$ {\it indecomposable} if $S$ is not
decomposable. \m
 
{\bf Example.}  It is possible for a linear IFS to have infinitely
many indecomposable eigensets corresponding to the same eigenvalue.
Consider $F = \{ \R^2; L_1, L_2 \}$ where
$$L_1 = \begin{pmatrix} 0 & -1 \\ 1 & 0 \end{pmatrix}, 
\qquad \qquad \qquad L_2 = \begin{pmatrix} 1 & 0 \\ 0 & 0.5 \end{pmatrix}.  
$$
Let $$S(r_1,r_2) = \{\, (\pm r_1, \pm r_2/2^k),\,  (\pm r_1, \mp r_2/2^k), \,
(\pm r_2/2^k, \pm r_1,), \, (\pm r_2/2^k, \mp r_1,) \, : \, k\geq 0\}.$$
It is easily verified that,
for any $r_1 \geq r_2 > 0$, the set $S(r_1,r_2)$ is an eigenset 
corresponding to eigenvalue $1$.  In addition, the unit square
with vertices $(1, 1),(1, -1),(-1, 1),(-1, -1) $ is also an eigenset 
corresponding
to eigenvalue $1$.\B

The proof of the following lemma is straightforward.
 A set $B \subset \R^n$ is called {\it star shaped} if $\lambda \, x
\in B$ for all for all $x\in B$ and all $0 \leq \lambda \leq 1$. 

\begin{lemma} \label{convex} 
\begin{enumerate}
\item If $\{ A_k\}$ is a sequence of centrally 
symmetric, convex, compact sets and $A$ is a compact set such that
$\lim_{k \rightarrow \infty} A_k = A$, then $A$ is 
also centrally symmetric and convex.
\item If $F$ is a compact, linear IFS, $B$ a centrally
symmetric, convex, compact set and $A = \lim_{k \rightarrow \infty}
F^k(B)$,  then $A$ is a centrally symmetric, star-shaped, compact set.
\end{enumerate}
\end{lemma}
 
\begin{lemma} \label{inclusion} If $F$ is an compact, irreducible, 
linear IFS with $\rho(F) = 1$, then there exists a compact, centrally
symmetric, convex body $A$ such that $F(A) \subseteq A$.
\end{lemma}

\proof Since, for each $k\geq 2$, we have $\rho ( (1-\frac{1}{k})F ) =
1-\frac{1}{k} < 1$, Theorem~\ref{linearCor} implies that there
is a compact, centrally symmetric, convex body $A_k$ such that
$$\left (1-\frac{1}{k} \right ) F(A_k) \subseteq int (A_k).$$ Since
$F$ is linear and the above inclusion is satisfied for $A_k$, it is
also satisfied for $\alpha \, A_k$ for any $\alpha >0$.  So, without
loss of generality, it can be assumed that $\max \{ \|x\| \, : x \in
A_k \} = 1$ for all $k\geq 2$.  Since the sequence of sets $\{A_k\}$
is bounded in $\mathbb H(\R^n)$, this sequence has an accumulation
point, a compact set $A$. Therefore, there is a subsequence $\{ A_{k_i} \}$
such that $\lim_{i\rightarrow \infty} \, A_{k_i} = A$ with respect to
the Hausdorff metric. Since
$$\left (1-\frac{1}{k_i} \right ) F(A_{k_i}) \subseteq int
(A_{k_i}),$$ it is the case that $\left (1-\frac{1}{k_i} \right )
f(A_{k_i}) \subseteq int (A_{k_i})$ for all $f \in F$. From this is is
straightforward to show that $f(A) \subseteq A$ for all $f\in F$ and
hence that $F(A) \subseteq A$.  Moreover, by Lemma~\ref{convex}, since
the $A_{k_i}$ are centrally symmetric and convex, so is $A$.  Notice
also that $A$ is a convex body, i.e., has nonempty interior; otherwise
$A$ spans a subspace $E\subset \R^n$ with $\dim \, E < n$ and  $F(A)
\subseteq A$ implies $F(E)
\subseteq E$, contradicting that $F$ is irreducible.  \qed \m 

The {\it affine span} aff($B$) of a set $B$ is the smallest affine
subspace of $\R^n$ containing $B$.  Call a set $B\subset \R^n$ {\it
full dimensional} if dim(aff$(B)) = n$.  Given an affine
IFS $F = \left( \R^n;
f_i,\, i\in I \right)$ let
$$F_\lambda = \left \{ \R^n ;\,\frac{1}{\lambda} f_i,\, i \in I \right \}.$$ 

\begin{lemma} \label{full} If an irreducible, affine IFS $F$ has
an eigenset $X$, then $X$ must be full dimensional.
\end{lemma}

\proof Suppose that $F(X) = \lambda X$, i.e. $F_{\lambda} (X) = X$.
For $x\in X$, let $g$ be a translation by $-x$.  For the IFS $F$, let
$F_g = \{ \R^n ; \, gfg^{-1}, f \in F_{\lambda} \}$.  If $Y= g(X)$,
then $0 \in Y$ and $F_g(Y) = Y$.  In particular, $Y$ is full
dimensional if and only if $X$ is full dimensional, and the affine
span of $Y$ equals the ordinary (linear) span $E = span(Y)$ of $Y$. 
Moreover, the linear parts of the affine maps in $F_g$ are just 
scalar multiples of the linear parts of the affine maps in $F$.
Therefore $F_g$ is irreducible if and only if $F$ is irreducible.

Let $f(x) = Lx +a$ be an arbitrary affine map in $F_g$. From
$F_g(Y) \subset Y \subset E$ it follows that $L(Y) + a = f(Y) \subset E$. 
Since $0\in Y$, also $a = L(0) + a = f(0) \in Y \subset E$.  Therefore
$L(Y) \subset -a + E = E$.  Since  $E = span(Y)$, also $L(E) \subset E$.
Because this is so for all $f \in F_g$, the subspace $E$
is invariant under all linear parts of maps in $F_g$.  
Because $F_g$ is irreducible, $dim(E) = n$.  Therefore
$Y$, and hence $X$, must be full dimensional.   
\qed \m 

\begin{lemma} \label{bounded} If 
$F = \left \{ \R^n ;\, L_i,\, i \in I \right \}$ is 
a bounded linear IFS, then there is an $\alpha >0$ such that $\alpha
\, F = \left \{ \R^n ;\, \alpha \, L_i,\, i \in I \right \}$ is
contractive.
\end{lemma}

\proof By the boundedness of $F$ there is an $R$ such that, for any
$L\in F$, $\frac{\|Lx \|}{\|x\|} \leq \|L\| \leq R$ for all $x\in
\R^n$.  Therefore, if $D_r$ denotes a disk of radius $r$ centered at
the origin, then $F(D_1) \subseteq D_R$. Hence $\frac{1}{2R} F(D_1)
\subset int(D_1)$.  By Theorem~\ref{linearCor} the IFS $\frac{1}{2R} \,F$ is
contractive.  \qed \B

\noindent {\bf Proof of Theorem~\ref{T-eigen}}: Given $F =
\left( \R^n; L_i,\, i\in I \right)$, consider the family
$\{F_\lambda\}$ of IFS's for $\lambda >0$.  Recall that $F_\lambda =
\left \{ \R^n ;\,\frac{1}{\lambda} f_i,\, i \in I \right \}$.

It is first proved that $F$ has no eigenvalue $\lambda > \rho(F)$.  By
way of contradiction assume that $\lambda > \rho(F)$, which implies
that $\rho(F_{\lambda}) < 1$.  According to Theorem~\ref{linearCor}
the IFS $F_{\lambda}$ is contractive.  By Corollary~\ref{unique} the
only invariant set of $F_{\lambda}$ is $\{0\}$, which means that the
only solution to the eigen-equation $F(X) = \lambda \, X$ is $X =
\{0\}$.  But by definition, $\{0\}$ is not an eigenset.  \m

The proof that $F$ has no eigenvalue $\lambda < \rho(F)$ is 
postponed because the more general affine version is provided
in the proof of Theorem~\ref{affThm} in Section~\ref{affSec}. \m

We now show that $\rho(F)$ is an eigenvalue of $F$.  Again let
 $F_{\lambda} = \frac{1}{\lambda}\, F$, so that $\rho(F_{\lambda})=
 1$. With $A$ as in the statement of Lemma~\ref{inclusion}, 
consider the nested intersection
$$S = \bigcap_{k\geq 0} F_{\lambda}^{k}(A) = \lim_{k \rightarrow
\infty} F_{\lambda}^k(A).$$ That $S$ is compact, centrally symmetric,
and star-shaped follows from Lemma~\ref{convex}.  Also
$$F_{\lambda}(S) = F_{\lambda} \left( \bigcap_{k\geq 0}
F_{\lambda}^{k}(A) \right ) = \bigcap_{k\geq 1} F_{\lambda}^{k}(A) =
S,$$ the last equality because $A \supseteq F_{\lambda}(A) \supseteq
F_{\lambda}^{(2)}(A) \supseteq \cdots$. From $F_{\lambda}(S)= S$ it
follows that $F(S) = \lambda \, S$.  \m

It remains to show that $S$ contains a non-zero vector.
 Since $A$ is a convex body and determined only up to scalar multiple,
 there is no loss of generality in assuming that $A$ contains a ball
 $B$ of radius $1$ centered at the origin. Then
$$\sup \, \{\, \|L_{\sigma}(x)\| \, : \, \sigma \in \Omega_k, \, x \in
B \} = \hat{\rho}_k(F_{\lambda}) \geq (\rho(F_{\lambda}))^k = 1.$$ So
there is a point $a_k \in F_{\lambda}^{k}(A)$ such that $\|a_k\| \geq
1$.  If $a$ is an accumulation point of $\{a_k\}$, then $\|a\| \geq
1$, and there is a subsequence $\{a_{k_i}\}$ of $\{a_{k}\}$ such that
$$\lim_{i \rightarrow \infty} a_{k_i} = a.$$ Since the sets
$F_{\lambda}^{(k_i)}(A)$ are closed and nested, it must be the case
that $a\in F_{\lambda}^{(k_i)}(A)$ for all $i$.  Therefore $a \in
S$. \m

That $S$ is full dimensional follows from Lemma~\ref{full}.
\qed 

\section{Theorems of Dranisnikov-Konyagin-Protasov and of Barabanov}
\label{Cor}

Important results of Dranisnikov-Konyagin-Protasov and of Barabanov on
the joint spectral radius turn out to be almost immediate corollaries
of Theorem~\ref{T-eigen}.  The first result is attributed to
Dranisnikov and Konyagin by Protasov, who provided a proof 
in \cite{P}.  Barabanov's theorem appeared originally in \cite{B}.  \m

\begin{cor} [Dranisnikov-Konyagin-Protasov] 
If $F = \left( \R^n; \, L_i,\, i \in I \right)$ is a compact,
irreducible, linear IFS with joint spectral radius $\rho := \rho(F)$,
then there exists a centrally symmetric convex body $K$ such that
$$conv \, F(K) = \rho K.$$ 
\end{cor}

\proof According to Theorem~\ref{T-eigen} there is a centrally
symmetric, full dimensional eigenset $S$ such that $F(S) = \rho \, S$.
If $K = conv (S)$, then $K$ is also centrally symmetric and
$$conv \, F(K) = conv \, F( conv \, S ) = conv \, F(S) = conv \, (\rho \, S) =
\rho \, conv \, S = \rho \, K.$$
The second equality is routine to check. Since $S$ is full dimensional,
$K$ is a convex body, i.e., has nonempty interior.
\qed \m 

The original form of the Barabanov theorem is as follows:  

\begin{theorem}[Barabanov] If a set $F$ of linear maps on
$\R^n$ is compact and irreducible, then there exists a vector norm
$\|\cdot \|_B$ such that
$$\begin{aligned} 
\mbox{for all $x$ and all $L\in F$} \qquad & 
\|Lx\|_B \leq \rho(F)\,\|x\|_B, \\
\mbox{for any $x\in \R^n$ there exists an $L\in F$ such that} \qquad & 
\|Lx\|_B = \rho(F)\,\|x\|_B.
\end{aligned}$$
\end{theorem}

Such a norm is called a {\it Barabanov norm}. The first property says
that $F$ is {\it extremal}, meaning that 
\begin{equation} \label{inequality} \|L\|_B \leq \rho(F) \end{equation} 
for all $L \in F$.  It is extremal
in the following sense. By property (4) of the joint spectral 
radius in Section~\ref{contractiveIFS}, 
$$\sup_{L\in F} \| L \| \geq \rho(F)$$ for any matrix norm.
Therefore, the joint spectral radius $\rho(F)$ can be
characterized as the infimum over all possible matrix norms of the
largest norm of linear maps in $F$.  Since $F$ is assumed compact, the
inequality~(\ref{inequality}) cannot be strict for all $L\in F$.  Hence
there exists an $L\in F$ whose Barabanov norm achieves the upper bound
$\rho(F)$.  Furthermore, the second property in the statement of
Barabanov's Theorem says that, for any $x\in \R^n$, there is such an
$L$ achieving a value equal to the joint spectral radius at the point
$x$. See \cite{W} for more on extremal norms. \m  

In view of Lemma~\ref{Minkowski}, Barabanov's theorem can
be restated in the following equivalent geometric form. Here
$\partial$ denotes the boundary.

\begin{cor} If $F$ is a compact, irreducible, linear IFS with joint
spectral radius $\rho := \rho(F)$, then there
exists a centrally symmetric convex body $K$ such that
$$F(K) \subseteq \rho K,$$ and, for any $x \in \partial K$, there
is an $L\in F$ such that $Lx \in \partial (\rho \, K)$.
\end{cor}

\proof Let $F^{t} = \left( \R^n; \, L^{t}_{i},\, i \in I\right)$, where
$L^{t}$ denotes the adjoint (transpose matrix) of $L$.  For a compact
set $Y$, the {\it dual} of $Y$ (sometimes called the polar) is the set
$$Y^* = \{ z \in \R^n \, : \langle y,z \rangle \leq 1 \hbox{ for all }
y \in Y \}.$$  The 
first two of the following
properties are easily proved for any compact set $B$.
\begin{enumerate}
\item $B^*$ is convex.
\item If $B$ is centrally symmetric, then so is $B^*$.
\item If $L$ is linear and $L^t(S) \subseteq  S$, then 
$L(S^*) \subseteq  S^*$.
\end{enumerate} 
To prove the third property above, assume that $L^t(S) \subseteq  S$.
and let $x\in S^*$.  Then 
$$\begin{aligned}
 x \in S^* & \Rightarrow \langle x,y \rangle \leq 1 \; 
\text{  for all  } y\in S \\ &
 \Rightarrow \langle x,L^ty \rangle \leq 1 \; \text{  for all  } y\in S \\
 & \Rightarrow \langle Lx,y \rangle \leq 1 \; \text{  for all  } y\in S \\
 & \Rightarrow Lx \in S^*
\end{aligned}$$

Since $F$ is a compact, irreducible, linear IFS, so is $F^t$.  Let $S$
be a centrally symmetric eigenset for $F^t$ as guaranteed by
Theorem~\ref{T-eigen}.  By properties 1 and 2 above, $S^*$ is a
centrally symmetric convex body.  From the eigen-equation $F^t(S) =
\rho \, S$, it follows that $\frac{1}{\rho} \, L^t(S) \subseteq S$ for all
$L\in F$.  From property 3 above it follows that 
$\frac{1}{\rho} \, F(S^*) \subseteq S^*$ or $ F(S^*) \subseteq \rho \, S^*$.
Setting $K = S^*$ yields $$F (K) \subseteq \rho \, K.$$

Concerning the second statement of the corollary, assume
that $x \in \partial K = \partial S^*$. Then $\langle x,y \rangle \leq
1$ for all $y\in S$ and $\langle x,y \rangle = 1$ for some $y\in
S$. Since $F(S) = \rho \, S$, the last equality implies that 
there is an $L \in F$ such that $\langle
\frac{1}{\rho} \, L x, z \rangle = \langle
x,\frac{1}{\rho} \, L^t z \rangle = 1$ for some $z\in S$. Now we have
$\langle \frac{1}{\rho} \, L x, y \rangle \leq 1$ for all $y\in S$ and
$\langle \frac{1}{\rho} \, L x, z \rangle = 1$ for some $z\in S$.
Therefore,  $\frac{1}{\rho} \, L x \in \partial S^* = \partial K$ or 
$L x \in \rho(\partial K) = \partial( \rho \, K)$.
 \qed

\section{The Eigenvalue Problem for an Affine IFS} \label{affSec}

For an affine IFS $F$, there is no theorem analogous to
Theorem~\ref{T-eigen}.  More specifically, there are examples where
$\rho(F)$ is an eigenvalue of $F$ and examples where $\rho(F)$ is not
an eigenvalue of $F$.  For an example where $\rho(F)$ is an
eigenvalue, let 

$$F_1 = \{ \R^2 ; \, f\}, \qquad  f(x) = Lx + (1,0), \qquad
L = \begin{pmatrix} 0 & -1 \\ 1 & 0
\end{pmatrix}.$$  Note that 
$L$, a $90^o$ degree rotation about the origin, is irreducible 
and $\rho(F_1) = 1$. If $S$ is the unit square with
vertices $(0,0), (1,0),(0,1),(1,1)$, then $F_1(S) = S$.  Therefore
$\rho(F_1) = 1$ is an eigenvalue of $F_1$.  On the other hand let 
$$F_2 = \{ \R ; \, f\}, \qquad \qquad f(x) = x + 1.$$ In this case
$\rho(F_2) = 1$, but it is clear that there exists no compact set $X$
such that $F(X) = X$.  For the affine case, Theorem~\ref{affThm}, as
stated in the introduction, does
holds.  The proof is as follows. \m

\noindent {\bf Proof of Theorem~\ref{affThm}}: 
If $\lambda > \rho(F)$, then $\rho(F_{\lambda}) < 1$.  According to
Theorem~\ref{aff}, the IFS $F_{\lambda}$ has an attractor $A$
so that $F_{\lambda}(A) = A$.  Since at least one function
in  $F_{\lambda}$ is not linear, $A \neq \{0\}$. Since $F_{\lambda}(A) = A$,
also $F(A) = \lambda \, A$.  Therefore $\lambda$ is an eigenvalue
of $F$. \m 

Concerning the second statement in the theorem assume, by way of
contradiction, that such an eigenvalue $\lambda < \rho(F)$ exists,
with corresponding eigenset $S$.  Then $F_{\lambda}(S) = S$ and
$\rho(F_{\lambda}) > 1$.  According to Lemma~\ref{full}, since $F$ is
assumed irreducible, the eigenset $S$ is full dimensional.  Exactly as
in the proof of Lemma~\ref{full}, using conjugation by a translation,
there is an affine IFS $F'$ and a nonempty compact set $S'$ such that
\begin{enumerate}
\item $F'(S') = S'$,
\item $0 \in int(conv(S'))$, 
\item The set $\mathbb L_{F'}$ of linear parts of the functions in $F'$
is equal to the set $\mathbb L_{F_{\lambda}}$ of linear parts of the functions
in $F_{\lambda}$, 
\item $\rho(F') = \rho(F_{\lambda}) > 1$,
\item $F'$ is irreducible
\end{enumerate}
In item 2 above, $int(conv(S'))$ denotes the interior of the convex
hull of $S'$.  If $K = conv(S')$ and $f(x) = Lx +a$ is an arbitrary
affine function such that $f(S') \subseteq S'$, then $$f(K) \subseteq
K.$$ This follows from the fact that $f(S') \subseteq S'$ as follows.
If $z \in K$, then $z = \alpha \, x + (1-\alpha)\, y$ where $0 \leq
\alpha \leq 1$ and $x,y \in S'$.  Therefore
$$\begin{aligned} f(z) & = \alpha \, Lx + (1-\alpha) \, Ly + a =
\alpha ( Lx +a) + (1-\alpha)( Ly + a) \\ & = \alpha \, f(x) +
(1-\alpha) \, f(y) \in conv(f(S')) \subset conv(S') = K. \end{aligned}$$ 

Let $r > 0$ be the largest radius of a ball centered at the origin and
contained in $K$ and $R$ the smallest radius of a ball centered at the
origin and containing $K$.  Let $x \in K$ such that $0 < \| x\| \leq r$.
If $f(x) = Lx + a$ is any affine function such that $f(S') \subseteq
S'$, then we claim that $\|Lx \| \leq R+r$.  To prove this, first note
that $-x \in K$. From
$f(K) \subseteq K$ it follows that
$$\begin{aligned} \| Lx + a \| & = \| f(x) \| \leq R \\
 \| -Lx + a \| & =  \| L(-x) + a \| = \| f(-x) \| \leq R \\
 \|2a\| & = \| (Lx + a) +  (-Lx + a) \| \leq  \| Lx + a \| + \| L(-x) + a \|
 \leq  2R \\ 
\|Lx \| & = \| f(x) - a \| \leq  \| f(x)\| + \|a \| \leq R+r.
 \end{aligned}$$ 

From the definition of the joint spectral radius, $\rho(F') > 1$
implies that there is an $\epsilon > 0$ such that $(\hat
{\rho}_k(F_{\lambda}))^{1/k} > 1+\epsilon$ for infinitely many values
of $k$.  This, in turn, implies that, for each such $k$, there is an
affine map $f_k \in \{ f_{\sigma} \, : \, \sigma \in \Omega_k \}$ and
its linear part $L_k \in \{ L_{\sigma} \, : \, \sigma \in \Omega_k \}$
such that $\|L_k \| \geq (1+\epsilon)^k$.  Choose $k = k_0$
sufficiently large that $\|L_k \| \geq (1+\epsilon)^{k_0} >
\frac{R+r}{r}$.  Then there is a $y \in K'$ with $\|y \| = r$ such
that $\| L_{k_0} y \| > r \, \frac{R+r}{r} = R+r$. Since $L_{k_0}$ is
the linear part of an affine function $f$ with the property $f(S')
\subseteq S'$ (property 1 above), this is a contradiction to what
was proved in the previous paragraph. \qed

\end{document}